\RequirePackage[hyphens]{url}
\documentclass[a4paper,USenglish]{lipics-v2021}
\usepackage[hideChapterNumber,refchapterCite,noCCBy]{craigstyle}

\usepackage{amsmath, amsthm, amssymb}
\usepackage[only,llbracket,rrbracket]{stmaryrd}

\usepackage{pifont}%
\newcommand{\cmark}{\ding{51}}%
\newcommand{\xmark}{\ding{55}}%
 
\usepackage{tabularx}
\usepackage{multirow}

\usepackage{tikz-cd}

\newtheorem{mydefinition}[theorem]{Definition}
\theoremstyle{definition}
\newtheorem{question}[theorem]{Question}
\newtheorem*{notation}{Notation}
\newtheorem*{convention}{Convention}

\newcommand{\liffdef}{\stackrel{\mathrm{def}}{\iff}}
\newcommand{\To}{\Rightarrow}

\newcommand{\mbf}[1]{\mathbf{#1}}

\newcommand{\mc}[1]{\mathcal{#1}}

\newcommand{\op}{{\mathrm{op}}}
\newcommand{\sem}[1]{\llbracket{#1}\rrbracket}

\newcommand{\isdef}{\stackrel{\mathrm{def}}{=}}

\newcommand{\into}{\hookrightarrow}

\newcommand{\bH}{\mathbf{H}}
\newcommand{\bE}{\mathbf{E}}
\DeclareMathOperator*{\Con}{\mathsf{Con}}
\newcommand{\rE}{\mathsf{E}}
\newcommand{\rA}{\mathsf{A}}

\DeclareMathOperator{\KCon}{\mathsf{KCon}}

\newcommand{\p}{\overline{p}}
\newcommand{\q}{\overline{q}}
\renewcommand{\r}{\overline{r}}

\newcommand{\x}{\overline{x}}
\newcommand{\y}{\overline{y}}
\renewcommand{\a}{\overline{a}}
\renewcommand{\b}{\overline{b}}
\newcommand{\Luklog}{\textbf{\L}}
\newcommand{\ThFO}{\mathrm{Th}_{\mathsf{FO}}}
\DeclareMathOperator{\mand}{\,\mathsf{and}\,}

\DeclareMathOperator{\Sub}{\mathsf{Sub}}
\renewcommand{\epsilon}{\varepsilon}

\newcommand*\cocolon{%
        \nobreak
        \mskip6mu plus1mu
        \mathpunct{}%
        \nonscript
        \mkern-\thinmuskip
        {:}%
        \mskip2mu
        \relax
}

\usepackage{tcolorbox}

\title{\vspace{26pt}Uniform Interpolation}
\author{Sam van Gool}{Université Paris-Saclay, ENS Paris-Saclay, CNRS, LMF, France}{}{}{}

\titlerunning{Uniform Interpolation}

\begin{document}

\maketitle

\begin{abstract}
   Uniform interpolation is a strengthening of interpolation that holds for certain propositional logics. The starting point of this chapter is a theorem of A. Pitts, which shows that uniform interpolation holds for intuitionistic propositional logic. We outline how this theorem may be proved semantically via the definability of bisimulation quantifiers, and how it generalizes to an open mapping theorem between Esakia spaces. We also discuss connections between uniform interpolation and research in categorical logic, algebra, and model theory.%
\end{abstract}

\tableofcontents

\section{Introduction}
Let $\psi(\p,\q)$ and $\theta(\q,\r)$ be propositional formulas and suppose that $\psi \to \theta$ is a tautology of classical logic.
Define the formula
\index{uniform interpolant!in classical propositional logic}
\begin{equation}\label{eq:booleanui}
    \phi(\q) \isdef \bigvee \left\{ \psi(\p/\b,\q) \ | \ \b \in \{\bot,\top\}^{\p} \right\},
\end{equation}
that is, $\phi$ is the disjunction of all possible variants of $\psi$ obtained by substituting some vector of truth values for the propositional variables $\p$.
As also explained in \refchapter{chapter:propositional}, Section~7, observe that $\phi$ is an interpolant for the entailment $\psi \vdash \theta$. Indeed, $\psi \vdash \phi$ follows from the definition of $\phi$, and $\phi \vdash \theta$ uses the assumption that $\psi \vdash \theta$.
Notice, moreover, that the interpolant $\phi$ defined here does not at all depend on the consequent formula $\theta$, but only on the antecedent formula $\phi$ and on the subset $\p$ of variables `to be eliminated' from $\psi$. In other words, $\phi$ is \emph{uniformly} an interpolant on the right of $\psi$: it interpolates any consequence $\theta'$ of $\psi$ that does not contain any of the variables in $\p$. The formula $\phi$ of~(\ref{eq:booleanui}) is therefore called a \emph{right uniform interpolant} for $\psi$ with respect to the set of eliminated variables $\p$.%
\footnote{We warn the reader of a possible source of confusion: \refchapter{chapter:propositional}, following part of the literature, calls this a right uniform $\q$-interpolant, thus referring in the name to the \emph{remaining} variables, rather than to the eliminated ones, as is done here, and also, for instance, in \refchapter{chapter:prooftheory}.}
\index{uniform interpolant!right}
Symmetrically, the expression
$\bigwedge \{\theta(\q/\b,\r) \ | \ \b \in \{\bot,\top\}^{\q}\}$
defines a \emph{left uniform interpolant} for $\theta$ with respect to $\p$ in classical logic. %
\index{uniform interpolant!left}

What other propositional logics have uniform interpolants? 
One may show that an entailment relation $\vdash$ has right uniform interpolants if it has interpolants  and is locally finite, by which we mean that there are finitely many equivalence classes of formulas in a finite set of variables. 
However, the local finiteness is by no means a \emph{necessary} condition for a propositional logic to have (right) uniform interpolants. Indeed, the following surprising theorem, which inspired much of the research surveyed in this chapter, is an instance of a non-locally-finite logic with uniform interpolation. %
\begin{theorem}[Pitts, \cite{Pit1992}] \label{thm:pitts}
    In intuitionistic logic, any propositional formula has both left and right uniform interpolants, with respect to any set of propositional variables. 
\end{theorem}
To give the general definition of uniform interpolation, we place ourselves in the context of a set $\mathcal{L}$ of formulas in a propositional logic, and we assume given a binary \emph{entailment relation} $\vdash$ between formulas in $\mathcal{L}$, which we will always assume is reflexive and transitive. 
\begin{mydefinition}\label{dfn:UI}
    Let $\vdash$ be an entailment relation, $\phi$ a formula, and $p$ a propositional variable. We call a formula $\rE_{p}(\phi)$ a \emph{right uniform interpolant} for $\phi$ with respect to $p$ if the following hold:
    \begin{enumerate}
        \item the formula $\rE_{p}(\phi)$ only contains variables that occur in $\phi$, but not $p$, 
        \item the entailment $\phi \vdash \rE_{p}(\phi)$ holds, and, 
        \item for any $p$-free formula $\theta$, if $\phi \vdash \theta$, then $\rE_{p}(\phi) \vdash \theta$.
    \end{enumerate}
Similarly, we call a formula $\rA_p(\phi)$ a \emph{left uniform interpolant} for $\phi$ with respect to $p$ if the following hold:
\begin{enumerate}
    \item the formula $\rA_{p}(\phi)$ only contains variables that occur in $\phi$, but not $p$, 
    \item the entailment $\rA_{p}(\phi) \vdash \phi$ holds, and, 
    \item for any $p$-free formula $\psi$, if $\psi \vdash \phi$, then $\psi \vdash \rA_{p}(\phi)$. 
\end{enumerate}
We say that the relation $\vdash$ \emph{has (left and right) uniform interpolants} if such formulas $\rE_p(\phi)$ and $\rA_p(\phi)$ exist for all formulas $\phi$ and propositional variables $p$.
\end{mydefinition}
We immediately make a few remarks about Definition~\ref{dfn:UI}. First, we have only given the definition with respect to a single variable, but the uniform interpolants with respect to sets of variables can be obtained by iterating the single-variable constructs. As a consequence, if an entailment relation has uniform interpolants with respect to single variables, then it has uniform interpolants with respect to any set of variables. 

Second, any formula $\phi$ has, up to equivalence, at most one right uniform interpolant, and at most one left uniform interpolant. Here, two formulas $\phi$ and $\phi'$ are \emph{equivalent} in an entailment relation $\vdash$ if both $\phi \vdash \phi'$ and $\phi' \vdash \phi$ hold. With a slight abuse of notation, we therefore often speak of `the' right uniform interpolant and `the' left uniform interpolant of $\phi$. 

Third, we have defined uniform interpolants with respect to an entailment relation $\vdash$, but there exists a subtly different definition in the literature. Indeed, let $\to$ be a binary symbol in the language under consideration. One may then also reasonably define%
\footnote{For instance, this is how uniform interpolation is defined in Definition~24 in Section~3.3 of \refchapter{chapter:prooftheory}.} \emph{implication-based uniform interpolants} by replacing, in Definition~\ref{dfn:UI}, all assertions of the form `$\alpha \vdash \beta$' by `$\top \vdash \alpha \to \beta$', where $\top$ is a fixed constant symbol in the language, suitably connected to the operation $\to$. Let us say that the entailment relation $\vdash$ has the \emph{deduction property} if, for all formulas $\alpha$ and $\beta$, $\alpha \vdash \beta$ is equivalent to $\top \vdash \alpha \to \beta$. If the entailment relation $\vdash$ does \emph{not} have this deduction property, then the notions of uniform interpolant in the sense of Definition~\ref{dfn:UI} and the implication-based uniform interpolant are not the same in general. Indeed, when $\vdash$ is the global consequence relation for the modal logic $\mathbf{K}$, uniform interpolants in the sense of Definition~\ref{dfn:UI} do not always exist, while implication-based uniform interpolants do exist~\cite[Sec.~6]{Ghi1995};  we come back to this in Remark~\ref{rem:K-UI} below. In light of this, the uniform interpolants of Definition~\ref{dfn:UI} are sometimes called \emph{deductive} uniform interpolants for emphasis.

\begin{convention}
In this chapter we adopt the convention that the phrase `uniform interpolant' always means `uniform interpolant in the sense of Definition~\ref{dfn:UI}', and we explicitly mention `implication-based uniform interpolant' when needed.
\end{convention}

Finally, note that the three conditions (1)--(3) defining the right and left uniform interpolants, respectively, correspond to rules for second order existential and universal quantifiers; this explains the suggestive notations $\rE_p(\phi)$ and $\rA_p(\phi)$. The connection between uniform interpolants and higher order quantifiers was central to Pitts' original paper~\cite[\S 4]{Pit1992}, and was further elaborated in the algebraic, categorical, and model-theoretic approaches to uniform interpolation, in particular by S. Ghilardi \& M. Zawadowski~\cite{GhiZaw2002}, as we will discuss in depth in Sections~\ref{sec:algebra} and \ref{sec:model-theory} below.
Note immediately, however, that the interpretation of second order quantifiers using the uniform interpolants $\rE_p$ and $\rA_p$  is not \emph{conservative}~\cite[Exa.~10]{Pit1992}; in other words, uniform interpolants provide \emph{a} possible interpretation of second order quantifiers, but moreover satisfy specific properties that generic second order quantifiers do not.

Proof methods for uniform interpolation can be divided, roughly, into two strands: one is syntactic and relies on the existence of a well-behaved sequent calculus for the entailment relation $\vdash$. This is how uniform interpolation was established for intuitionistic entailment in Pitts' original proof~\cite{Pit1992}.
Pitts' proof relies on a strongly terminating proof calculus for intuitionistic propositional logic, and has been the subject of intensive proof-theoretic investigations and extensions to other logics; see, for instance, Section~3.3 in \refchapter{chapter:prooftheory} and the references therein; we will also make some more remarks on this method in Section~\ref{sec:connections} of this chapter.
The other proof method, which will be our focus in this chapter, starts from a \emph{semantic} characterization of the entailment relation $\vdash$, typically in terms of a consequence relation on Kripke models. 
One then shows that \emph{bisimulation quantifiers} are definable in the propositional language, and deduces from this the existence of {implication-based} uniform interpolants. This approach to proving {implication-based} uniform interpolation was developed soon after~\cite{Pit1992} as an alternative proof for Theorem~\ref{thm:pitts}~\cite{GhiZaw1995}, and was also shown to apply to modal logics such as Gödel-Löb logic $\mathbf{GL}$~\cite[Cor.~2.13]{Shavrukov}, modal logic $\mathbf{K}$~\cite[Thm~6.1]{Ghi1995} and Grzegorczyk logic $\mathbf{S4Grz}$~\cite[Sec.~8]{Vis1996}.

The rest of this chapter is structured as follows. In Section~\ref{sec:bisimulation}, we outline the proof of uniform interpolation for intuitionistic logic via bisimulation quantifiers. We then make a link with adjoints between free algebras, and with Esakia duality, in Section~\ref{sec:algebra}. In Section~\ref{sec:model-theory}, we expose the relationship between uniform interpolants and quantifier elimination in first order theories of classes of algebras associated to propositional logics. In Section~\ref{sec:connections}, we point out some further connections, in particular to syntactic methods, and give a table, Table~\ref{tab:summary}, summarizing what is (not) known about uniform interpolation for some of the most commonly considered propositional logics.

Note that our focus in this chapter is on uniform interpolation for \emph{propositional} logics. In the classical setting, first order logic does not have uniform interpolants, see for instance Theorem~2.17 in \refchapter{chapter:firstorder}. There are substructural predicate logics with uniform interpolation~\cite{AliDerOno2014}. We do not know at present whether or not intuitionistic predicate logic has uniform interpolation.

\section{Uniform Interpolants via Bisimulation Quantifiers}
\label{sec:bisimulation}
In this section, we explain how to prove Theorem~\ref{thm:pitts} using the idea of  \emph{bisimulation quantifiers}. This proof of Theorem~\ref{thm:pitts} is essentially due to~\cite{GhiZaw1995}, which was inspired by a similar technique developed for $\mathbf{GL}$ in~\cite{Shavrukov}. While~\cite{GhiZaw1995} phrases the argument in terms of sheaf representations of finitely presented Heyting algebras, we here present the ideas in more elementary terms, following~\cite{Vis1996}; for more about the connection, see Subsection~\ref{subsec:further-bisim} below.

We focus on the case of intuitionistic propositional logic for our presentation here, but similar ideas apply for proving implication-based uniform interpolation in modal logics; also see Subsection~\ref{subsec:further-bisim} for more about this.

\subsection{Models and Bisimilarity}
For $\p$ a set of propositional variables, an (intuitionistic) \emph{Kripke model over $\p$}, or \emph{$\p$-model}, is a tuple $M = (W, \preceq, \pi)$, where $W$ is a set of \emph{nodes}, $\preceq$ is a reflexive and transitive binary relation on $W$, and $\pi$ is a function that assigns to every variable $p \in \p$ an upwards closed subset $\pi(p)$ of $W$; such a model is called \emph{finite} if its set of nodes is finite. A \emph{pointed $\p$-model} is a tuple $(M, w)$, where $M = (W, \preceq, \pi)$ is a $\p$-model and $w \in W$. For any pointed $\p$-model $(M,w)$ and any formula $\phi(\p)$, we write $M, w \Vdash \phi$ if $\phi$ is \emph{forced} at node $w$. This forcing relation $\Vdash$ is defined inductively, as usual, for any pointed $\p$-model $(M, w)$, by the following five clauses, where $p$ is any variable in $\p$ and $\phi_1, \phi_2$ are any formulas with variables in $\p$:
\begin{enumerate}
    \item $M, w \Vdash p \liffdef w \in \pi(p)$;
    \item $M, w \Vdash \phi_1 \wedge \phi_2 \liffdef M, w \Vdash \phi_1 \text{ and } M,w \Vdash \phi_2$;
    \item $M, w \Vdash \phi_1 \vee \phi_2 \liffdef M, w \Vdash \phi_1 \text{ or } M,w \Vdash \phi_2$;
    \item $M, w \Vdash \phi_1 \to \phi_2 \liffdef $ for any node $v$ in $M$ such that $w \preceq v$, if $M, v \Vdash \phi_1$ then $M, v \Vdash \phi_2$;
    \item $M, w \not\Vdash \bot$.
    \end{enumerate}

    Intuitionistic entailment is \emph{sound and complete} with respect to finite Kripke models. More precisely, if $\phi(\p)$ and $\psi(\p)$ are formulas, then the intuitionistic entailment $\phi \vdash \psi$ holds if, and only if, for any finite pointed $\p$-model $(M,w)$, if $M, w \Vdash \phi$ then $M, w \Vdash \psi$. 
    For a proof of this foundational equivalence, see, e.g., \cite[Thm.~2.43]{ChaZak1997}.
    In light of this equivalence, when $\phi(\p)$ is a formula and $\q$ is any set of variables with $\p \subseteq \q$, we define 
\[ \sem{\phi}_{\q} := \{ (M,w) \text{ a finite pointed $\q$-model} \ | \ (M, w) \Vdash \phi \}\]
and we note that soundness and completeness means that $\phi \vdash \psi$ is equivalent to $\sem{\phi}_{\q} \subseteq \sem{\psi}_{\q}$.

We will make use of the key semantic notion of \emph{bisimilarity}, a combinatorial method for determining when two nodes in Kripke models have the same `behavior'. In our setting, having the same behavior means that the two nodes force the same intuitionistic propositional formulas. In order to prove the existence of uniform interpolants with respect to a propositional variable $p$, we will in particular be interested in nodes that have the same behavior \emph{for a given set of variables}. The following notion of $\q$-bisimilarity captures this idea.

\begin{mydefinition}
Let $M = (W, \preceq, \pi)$ and $M' = (W', \preceq', \pi')$ be $\p$-models, and let $\q \subseteq \p$. A \emph{$\q$-bisimulation} between $M$ and $M'$ is a relation $Z \subseteq W \times W'$ satisfying the following three properties, for any $(w,w') \in Z$:
\begin{enumerate}
    \item for any $q \in \q$,  we have $w \in \pi(q)$ if, and only if, $w' \in \pi'(q)$;
    \item for any $v \in W$, if $w \preceq v$, then there exists $v' \in W'$ such that $(v,v') \in Z$ and $w' \preceq' v'$;
    \item for any $v' \in W'$, if $w' \preceq' v'$, then there exists $v \in W$ such that $(v,v') \in Z$ and $w \preceq v$.
\end{enumerate}
We say that a pair of nodes $(w,w') \in W \times W'$ is \emph{$\q$-bisimilar} if there exists a $\q$-bisimulation between $M$ and $M'$ which contains the pair $(w,w')$. In this case, we write $(M,w) \sim_{\q} (M',w')$.\footnote{We warn the reader that, in some sources, the relation that we denote by $\sim_{\q}$ would instead be denoted by $\sim_{\r}$, where $\r = \p \setminus \q$, thus  emphasizing the set of propositional variables that have been `forgotten', rather than the set of variables that have been `kept'.}
\end{mydefinition}
A fundamental fact is that $\q$-bisimilar nodes force exactly the same formulas in variables $\q$: %
\begin{lemma}\label{lem:bisim-form}
    Let $(M,w)$ and $(M',w')$ be pointed $\p$-models, let $\q \subseteq \p$, and suppose that $(M, w) \sim_{\q} (M',w')$. For any formula $\phi(\q)$, we have $(M, w) \Vdash \phi$ if, and only if, $(M', w') \Vdash \phi$.
\end{lemma}
\begin{proof}
    By induction on the construction of the formula $\phi$.
\end{proof}

\subsection{Bisimulation Quantifiers}
The \emph{existential bisimulation quantifier} is given by the following semantic construction. Let $p,\q$ be a set of variables and let $\mathcal{K}$ be a class of finite pointed $(p, \q)$-models. We define the class of finite pointed $(p,\q)$-models
\[ \mathcal{E}_p(\mathcal{K}) := \{ (M, w) \ | \ \text{there exists } (M', w') \text{ such that } (M, w) \sim_{\q} (M',w') \text{ and } (M',w') \in \mathcal{K} \}. \] 
Similarly, the \emph{universal bisimulation quantifier} is defined as the class of finite pointed $(p,\q)$-models
\[ \mathcal{A}_p(\mathcal{K}) := \{ (M, w) \ | \ \text{for all } (M', w') \text{ such that } (M, w) \sim_{\q} (M',w'), \text{ we have } (M', w') \in \mathcal{K} \}. \]

The following proposition shows that \emph{definability} of bisimulation quantifiers implies the existence of uniform interpolants. We give a simple proof in which we already assume the usual, non-uniform Craig interpolation theorem for intuitionistic propositional logic. It is possible to avoid the appeal to this theorem, as we will explain after the proof.
\begin{proposition}\label{prop:definability}
Let $\phi(p,\q)$, $\epsilon(\q)$, and $\alpha(\q)$ be formulas. 
\begin{enumerate}
    \item If  $\sem{\epsilon}_{p,\q} = \mathcal{E}_p(\sem{\phi}_{p,\q})$, then $\epsilon$ is a right uniform interpolant for $\phi$ with respect to $p$.
    \item If $\sem{\alpha}_{p,\q} = \mathcal{A}_p(\sem{\phi}_{p,\q})$, then $\alpha$ is a left uniform interpolant for $\phi$ with respect to $p$.
\end{enumerate}
\end{proposition}
\begin{proof}
    (1) Assume that $\sem{\epsilon}_{p,\q} = \mc{E}_p(\sem{\phi}_{p,\q})$. 
    Note that $\sem{\phi}_{p,\q} \subseteq \mc{E}_p(\sem{\phi}_{p,\q})$, since any model is $p$-bisimilar to itself. Thus, $\sem{\phi}_{p,q} \subseteq \sem{\epsilon}_{p,\q}$, and therefore $\phi \vdash \epsilon$ by completeness. Now let $\theta(\q,\r)$ be any formula such that $\phi \vdash \theta$. By the Craig interpolation theorem for intuitionistic logic, pick $\theta'(\q)$ such that $\phi \vdash \theta'$ and $\theta' \vdash \theta$. We will prove that $\epsilon \vdash \theta'$, from which $\epsilon \vdash \theta$ follows by transitivity of $\vdash$. To establish $\epsilon \vdash \theta'$, by completeness it suffices to show that $\sem{\epsilon}_{p,\q} \subseteq \sem{\theta'}_{p,\q}$. Let $(M,w) \in \sem{\epsilon}_{p,\q}$ be arbitrary. 
    Since $\sem{\epsilon}_{p,\q} = \mc{E}_p(\sem{\phi}_{p,\q})$, pick $(M',w') \sim_{\q} (M,w)$ such that $M',w' \Vdash \phi$. Since $\phi \vdash \theta'$, by soundness we have $M',w' \Vdash \theta'$. By Lemma~\ref{lem:bisim-form}, we have $M,w \Vdash \theta'$, as required. The proof of (2) is analogous.
\end{proof}
    The converse statements of (1) and (2) in Proposition~\ref{prop:definability} are also true, as we will prove in Corollary~\ref{cor:converse-definability} below. We also remark that the appeal to Craig interpolation theorem in the proof of Proposition~\ref{prop:definability} could be avoided, by modifying the statement in item (1) to: if, for any finite set of variables $\r$, $\sem{\epsilon}_{p,\q,\r} = \mathcal{E}_p(\sem{\phi}_{p,\q,\r})$, then $\epsilon$ is a right uniform interpolant for $\phi$ with respect to $p$. The proof of this modified statement is very similar to the proof of Proposition~\ref{prop:definability}(1), but does not require the Craig interpolation theorem. A similar remark holds for item (2). Using this modified version of Proposition~\ref{prop:definability}, one may slightly strengthen the proof outlined below and give a proof of the uniform interpolation theorem that does not already rely on the Craig interpolation theorem.

A more abstract point of view on bisimulation quantifiers $\mc{E}_p$ and $\mc{A}_p$ is that they are semantic realizations of left and right \emph{adjoints}, respectively, of certain homomorphisms between free Heyting algebras. Proposition~\ref{prop:definability} then says that if these adjoints are definable in the logic, then they give uniform interpolants. We will expand upon this idea in Section~\ref{sec:algebra} below.

\subsection{Definability via Bounded Bisimulations}
In light of Proposition~\ref{prop:definability}, we now seek to prove that the existential and universal bisimulation quantifiers are definable in intuitionistic propositional logic. To do so, the main technical tool is that of \emph{bounded bisimulation}~\cite{Vis1996} or, equivalently, the \emph{bisimulation games} of~\cite{GhiZaw1995}.
\begin{mydefinition}
    Let $M = (W, \preceq,\pi)$ and $M' = (W', \preceq', \pi')$ be $\q$-models. Let $Z_0, \dots, Z_n \subseteq W \times W'$ be a finite sequence of binary relations. We call $(Z_i)_{i=0}^n$ a \emph{bounded bisimulation} of depth $n$ if it satisfies the following three properties:    \begin{enumerate}
        \item for any $(w,w') \in \bigcup_{i=0}^n Z_i$ and $q \in \q$, $w \in \pi(q)$ if, and only if, $w' \in \pi'(q)$,
\end{enumerate}
and, for each $0 \leq i < n$ and $(w,w') \in Z_{i+1}$,
\begin{enumerate}
    \setcounter{enumi}{1}
        \item for any $v \in W$, if $w \preceq v$, then there exists $v' \in W'$ such that $(v,v') \in Z_i$ and $w' \preceq' v'$;
        \item for any $v' \in W'$, if $w' \preceq' v'$, then there exists $w \in W$ such that $(w,w') \in Z_i$ and $w \preceq v$.
    \end{enumerate}
    We say that $M,w$ and $M',w'$ are \emph{$n$-bisimilar} if there exists a bounded bisimulation of depth $n$ which relates $w$ to $w'$. In this case, we write $M, w \sim_{n,\q} M', w'$.
\end{mydefinition}
Note that if $Z$ is a bisimulation, then the constant sequence $Z, \dots, Z$ gives a bounded bisimulation of arbitrary depth. Thus, $M, w \sim_{\q} M',w'$ implies that, for every $n \geq 0$, $M, w \sim_{n,\q} M',w'$. The converse does not hold in general, but it does hold for finite models.\footnote{Indeed, we will see below that `for all $n \geq 0$, $M, w \sim_{n,\q} M',w'$' means precisely that $(M,w)$ and $(M',w')$ satisfy the same intuitionistic formulas. 
\cite[Cor.~35]{Pat1997} shows that for finite models, this is equivalent to bisimilarity, while an example in~\cite[Prop.~27]{Pat1997} shows that it is not equivalent in general. When logical equivalence implies bisimilarity, this is called a \emph{Hennessy-Milner} property; for more about such properties in the context of intuitionistic logic, see~\cite{GroPat2022} and the references therein.}

We also introduce the following convenient notation:
\begin{itemize}
    \item $M, w \preceq_{0,\q} M', w'$ means that, for every $q \in \q$, if $w \in \pi(q)$, then $w' \in \pi'(q)$; 
    \item for any $n \geq 1$, $M, w \preceq_{n,\q} M', w'$ means that, for every $v'$ in $M'$ such that $v' \succeq' w'$, there exists $v$ in $M$ with $v \succeq w$ and $M, v \sim_{n-1,\q} M', v'$.
\end{itemize}
With these notations, observe that, for any pointed $\q$-models $(M,w), (M',w')$ and any $n \geq 0$, we have  $M, w \sim_{n,\q} M', w'$ if, and only if, both $(M, w) \preceq_{n,\q} (M', w')$ and $M', w' \preceq_{n,\q} M, w$ hold.

For readers familiar with the bisimulation games of~\cite{GhiZaw1995, GhiZaw2002}, we note that pointed models $M, w$ and $M', w'$ are $n$-bisimilar precisely when Player II has a winning strategy for the first $n$ moves in the bisimulation game played on these models.

The interest of bounded bisimulation is its close relationship to definability of classes of pointed models. We formalize this in Propositions~\ref{prop:theories}~and~\ref{prop:bounded-form} below. Define the \emph{depth} of a propositional formula $\phi$ to be the maximum nesting depth of the implication connective $\to$. For any finite set of variables $\p$ and $n \geq 0$, write $F_{n,\p}$ for the set of intuitionistic equivalence classes of formulas in the variables $\p$ of depth $\leq n$. Note that, for any $n$ and $\p$, the set $F_{n,\p}$ is finite, as may be shown by a simple induction on $n$ using disjunctive normal forms. Fix, for each class $C \in F_{n,\p}$, a representing formula $\phi_C(\p)$ of depth $\leq n$. We will often confuse a class with its representing formula. In particular, when $\tau$ is a subset of $F_{n,\p}$, we will write $\bigwedge \tau$ to denote the formula $\bigwedge_{C \in \tau} \phi_C$. %

For any pointed $\q$-model $(M,w)$ and $\p \subseteq \q$, define 
\[ \mathrm{Th}_{n,\p}(M,w) := \{ C \in F_{n,\p} \ | \ M, w \Vdash \phi_C \},\] 
where we note that whether or not a class $C$ belongs to $\mathrm{Th}_{n,\p}(M,w)$ does not depend on the particular choice of representing formula $\phi_C$.
\begin{proposition}\label{prop:theories}
    Let  $(M,w)$ and $(M',w')$ be finite pointed $\q$-models. For any $n \geq 0$, we have $M, w \preceq_{n,\q} M', w'$ if, and only if,     $\mathrm{Th}_{n,\q}(M,w) \subseteq \mathrm{Th}_{n,\q}(M',w')$.
\end{proposition}
\begin{proof}
    See, e.g., \cite[Sec.~4]{Vis1996}, compare also \cite[Sec.~2]{GhiZaw1995}.
\end{proof}
Note that Proposition~\ref{prop:theories} and the preceding remarks in particular entail that for any $n, \q$, the relation $\sim_{n,\q}$ is an equivalence relation on the class of finite pointed $\q$-models. Moreover, this equivalence relation has only finitely many equivalence classes, because each $\sim_{n,\q}$-class uniquely determines a subset of $F_{n,\q}$, of which there are finitely many.

This idea lets us relate $\sim_{n,\q}$-bisimilarity to definability, as we will show now. In order to conveniently state the following proposition, if $\mathcal{K}$ is a class of finite pointed $\p$-models, then for $n \geq 0$ and $\q \subseteq \p$, we will say that the class $\mathcal{K}$ is \emph{$(n,\q)$-upward closed} if, for any $(M, w) \in \mc K$, and any finite pointed $\p$-model $(M',w')$, if $M,w \preceq_{n,\q} M',w'$, then $(M',w') \in \mc K$.
\begin{proposition}\label{prop:bounded-form}
    Let $\p$ be a finite sequence of variables, and let $\mc K$ be a class of finite pointed $\p$-models. For any $n \geq 0$ and $\q \subseteq \p$, the class $\mc K$ is $(n,\q)$-upward closed if, and only if, there exists a formula $\chi(\q)$ of depth $\leq n$ such that $\sem{\chi}_{\p} = \mc K$.
\end{proposition}
\begin{proof}
   The sufficiency is immediate from Proposition~\ref{prop:theories}. For the necessity, assume $\mc{K}$ is $(n,\q)$-upward closed, and consider the set $T_{\mc{K}} := \{ \mathrm{Th}_{n,\q}(M,w) \ | \ (M,w) \in \mc{K} \}$, which is a subset of $\mathcal{P}(F_{n,\q})$, and thus finite. 
    Define 
    \[ \chi := \bigvee_{\tau \in T_{\mc{K}}} \bigwedge \tau \ , \] 
    which is clearly a formula of depth $\leq n$ that only uses variables in $\q$. We show that $\sem{\chi}_{\p} = \mc{K}$. 
    First, if $(M, w) \in \mc{K}$, then clearly $M, w \Vdash \bigwedge \mathrm{Th}_{n,\q}(M,w)$, so $M, w \Vdash \chi$. Conversely, suppose that $M',w' \Vdash \chi$. Pick $(M,w) \in \mc{K}$ such that $M',w' \Vdash \bigwedge \mathrm{Th}_{n,\q}(M,w)$. This means that $\mathrm{Th}_{n,\q}(M, w) \subseteq \mathrm{Th}_{n,\q}(M', w')$. By Proposition~\ref{prop:theories}, it follows that $M, w \preceq_{n,\q} M', w'$. Since $\mc{K}$ is $(n,\q)$-upward closed, we conclude that $(M', w') \in \mc{K}$.
\end{proof}
Combining Proposition~\ref{prop:bounded-form} with Proposition~\ref{prop:definability}, we see that, in order to establish the existence of uniform interpolants, one needs to show that, for any $\phi(p,\q)$, the sets $\mc{E}_p(\sem{\phi}_{p,\q})$ and $\mc{A}_p(\sem{\phi}_{p,\q})$ are $(n,\q)$-upward closed, for some $n \geq 0$. This is the content of the final lemma, the proof of which is a combinatorial argument on Kripke models.
\begin{lemma}[Expansion Lemma]\label{lem:combi}
For any $k \geq 0$, there exists $n$ such that, whenever $\phi(p,\q)$ is a formula of depth $\leq k$, the sets  $\mc{E}_p(\sem{\phi}_{p,\q})$ and $\mc{A}_p(\sem{\phi}_{p,\q})$ are $(n,\q)$-upward closed.
\end{lemma}
We emphasize that Lemma~\ref{lem:combi} contains the combinatorial core of the uniform interpolation theorem, when viewed semantically. This is the main part of the argument that is highly specific to the class of models (and thus to the logic) under consideration, and often requires an intricate combinatorial argument; we refer to \cite[Lem.~4.2]{GhiZaw1995}, \cite[Lem.~5.1]{Vis1996}, or \cite[Lem.~10]{GooReg2017} for various renderings of such an argument in the case of intuitionistic logic. Note that, in all of these semantic arguments, the number $n$ in general needs to be taken very large compared to $k$.  The name `Expansion Lemma' comes from \cite[Lem.~2.10]{Shavrukov}, where it is used for the analogous combinatorial result that is required for proving uniform interpolation in Gödel-Löb logic $\mathbf{GL}$. 

Given the above results, we can now also see that any right uniform interpolant for $\phi$ with respect to $p$ must in fact have $\mathcal{E}_p(\sem{\phi}_{p, \q})$ as its semantics; that is, we deduce the converse of Proposition~\ref{prop:definability}.
\begin{corollary}\label{cor:converse-definability}
    Let $\phi(p,\q)$ be a formula.
    \begin{enumerate}
    \item If $\epsilon(\q)$ is a right uniform interpolant for $\phi$ with respect to $p$, then $\sem{\epsilon}_{p,\q}  = \mathcal{E}_p(\sem{\phi}_{p,\q})$.
    \item If $\alpha(\q)$ is a left uniform interpolant for $\phi$ with respect to $p$, then $\sem{\alpha}_{p,\q} = \mathcal{A}_p(\sem{\phi}_{p,\q})$.
\end{enumerate}
\end{corollary}
\begin{proof}
    (1) Assume that $\epsilon$ is a right uniform interpolant for $\phi$ with respect to $p$. By Proposition~\ref{prop:bounded-form} and Lemma~\ref{lem:combi}, pick a formula $\chi(\q)$ such that $\sem{\chi}_{p,\q} = \mathcal{E}_p(\sem{\phi}_{p,\q})$. By Proposition~\ref{prop:definability}, $\chi$ is a right uniform interpolant for $\phi$. It follows from the second remark following Definition~\ref{dfn:UI} that $\chi$ and $\epsilon$ are equivalent. In particular, $\sem{\epsilon}_{p,\q} = \sem{\chi}_{p,\q} = \mathcal{E}_p(\sem{\phi}_{p,\q})$. The proof of (2) is analogous.
\end{proof}
Corollary~\ref{cor:converse-definability} depends on the choice of the notion of `model'. It can happen that a logic with uniform interpolants is complete with respect to a class of models, but the bisimulation quantifiers fail to be definable when restricted to this class; an example in the case of $\mbf{GL}$ is given in \cite[Exa.~2, p. 108]{Ago2007}.

\subsection{Extensions of the Semantic Method to Other Logics}\label{subsec:further-bisim}
The proof of uniform interpolation for intuitionistic logic outlined above suggests a general proof strategy for uniform interpolation via definability of bisimulation quantifiers for logics admitting a Kripke semantics; note in particular that the argument given in Proposition~\ref{prop:definability} does not use any specific properties of intuitionistic logic beyond completeness and usual Craig interpolation (and the latter could have been avoided with some more care). The same proof method is for instance used in \cite{Vis1996} to obtain uniform interpolants in Grzegorczyk logic $\mathbf{S4Grz}$, and it is also shown there that implication-based uniform interpolants for $\mathbf{K}$~\cite[Sec.~6]{Ghi1995} and uniform interpolants for $\mathbf{GL}$~\cite{Shavrukov} can be obtained in this way.  %

The article~\cite{GhiZaw1995}, which preceded~\cite{Vis1996}, gave a very similar argument for uniform interpolation in intuitionistic logic, phrased in a more category-theoretic language. The basic idea in~\cite{GhiZaw1995} (later much expanded in~\cite{GhiZaw2002}) is that a formula $\phi(\q)$ can be represented as a \emph{sheaf} (i.e., a particular kind of functor), $F_{\q}(\phi)$, from the category of finite pointed posets with bounded morphisms to the category of sets. This sheaf $F_{\q}(\phi)$ associates to any finite pointed poset $(W, \preceq, w)$ the set of persistent valuations $\pi \colon \q \to \mathcal{P}(W)$ such that the model $((W, \preceq, \pi), w)$ forces $\phi$. In this setting, existential and universal bisimulation quantifiers with respect to a variable $p$ then correspond to certain {natural transformations} between these sheaves. The existence of uniform interpolants is obtained by showing that these natural transformations \emph{have an index}, a concept that is analogous to what we called `$(n,\q)$-upward closed' above. Inspired by this work, \cite{GooReg2017} gives a topological proof of Theorem~\ref{thm:pitts} via Esakia duality, see Subsection~\ref{subsec:topology} below. We will go into further detail about the work of \cite{GhiZaw1995,GhiZaw2002} in Section~\ref{sec:model-theory}. %

The semantic method also allows one to prove that uniform interpolation \emph{fails} to hold for certain logics. The first instance of this is~\cite{GhiZaw1995a}, where it is shown that uniform interpolants fail to exist for the modal logic $\mathbf{S4}$; \cite[Sec.~3]{Bil2007} shows that the same counterexample can be used to show that uniform interpolants do not exist for $\mathbf{K4}$. Uniform interpolation also fails for intuitionistic versions of these logics, $\mathbf{iK4}, \mathbf{iS4}$~\cite[Sec.~3.4]{Gie22}. %

\begin{remark}\label{rem:K-UI}
    The difference between implication-based and deductive uniform interpolants can be seen clearly in the case of the modal logic $\mathbf{K}$. An algorithm for computing implication-based uniform interpolants for the logic $\mathbf{K}$ is given in~\cite{Bil2007}, and is also outlined in Subsection~3.3 of \refchapter{chapter:prooftheory}. The fact that \emph{deductive} uniform interpolants do not always exist in $\mathbf{K}$ is essentially due to \cite[Thm.~37]{GhiLutWol2006long}, where a counterexample is given in the syntax of description logic. The failure is closely related to the failure of the deduction theorem for the global consequence relation, for which see, e.g.,~\cite[Sec.~1.5]{BRV2001}. We here give a slightly simpler example, due to G.~Metcalfe (personal communication); essentially the same example can be used to show that any strongly complete modal logic which has right uniform interpolants must be weakly transitive. Here, one may even weaken the assumption of having right uniform interpolants to being coherent in the sense of Theorem~\ref{thm:coherence} below~\cite[Thm.~4.1]{KowMet2018}.

    Consider the modal formula $\phi := (q \to p) \wedge (p \to \Box p) \wedge (p \to r)$. We claim that $\phi$ does not have a right uniform interpolant with respect to $p$ for the global consequence relation $\vdash_{\mbf K}$. The crucial observation towards this claim is that, for any modal formula $\theta(q,r)$,
    \begin{equation}\label{eq:K-example}
    \phi \vdash_{\mbf K} \theta \text{ if, and only if, there exists } N \geq 0 \text{ such that } \bigwedge_{k = 0}^N (q \to \Box^k r) \vdash_{\mbf K} \theta.
    \end{equation} 
    For the right-to-left direction of (\ref{eq:K-example}), just note that $\phi \vdash_{\mbf K} q \to \Box^k r$ for every $k \geq 0$ and use transitivity of $\vdash_{\mbf K}$. For the left-to-right direction of (\ref{eq:K-example}), we reason contrapositively, and assume that there does not exist $N$ such that $\bigwedge_{k = 0}^N (q \to \Box^k r) \vdash_{\mbf K} \theta$. Using the compactness and completeness theorems for modal logic $\mbf K$, there exists a pointed (modal) Kripke $(q,r)$-model $(M, w)$ such that $M,w \Vdash q \to \Box^k r$ for every $k \geq 0$, but $M,w \not\Vdash \theta$. Extend the valuation of this model to the variable $p$ by setting $M, w \Vdash p$ if, and only if, $M,w \Vdash \Box^k r$ for every $k \geq 0$. It it then straight forward to check that $M, w \Vdash \phi$, so that $\phi \not\vdash_{\mbf K} \theta$ by soundness.
    
    Now suppose that $\epsilon(q,r)$ would be a deductive right uniform interpolant for $\phi$ with respect to $p$. Since $\phi \vdash_{\mbf K} \epsilon$, by (\ref{eq:K-example}) we can pick $N \geq 0$ such that $\bigwedge_{k=0}^N (q \to \Box^k r) \vdash_{\mbf K} \epsilon$. Also, since $\phi \vdash_{\mbf K} q \to \Box^{N+1} r$, we must have $\epsilon \vdash_{\mbf K} q \to \Box^{N+1} r$. By transitivity of $\vdash_{\mbf K}$, we now obtain that $\bigwedge_{k=0}^N (q \to \Box^k r) \vdash_{\mbf K} q \to \Box^{N+1} r$, but this is false, as can be seen from a simple countermodel construction.\lipicsEnd
\end{remark}

Bisimulation quantifiers were also used to prove that the modal $\mu$-calculus has uniform interpolation~\cite[Cor.~3.8]{AgoHol2000}. The proof of this result relies on the characterization of the sets of transition systems that are definable in the $\mu$-calculus as being precisely those that are recognizable by so-called $\mu$-automata~\cite{JanWal1996}. The definability of bisimulation quantifiers is then obtained by showing that $\mu$-automata are closed under projections. This method is further used (\cite[Thm.~3.9]{AgoHol2000}) to prove that adding bisimulation quantifiers to the syntax of propositional dynamic logic yields a logic that is expressively equivalent to the full modal $\mu$-calculus. An alternative, proof-theoretic approach to uniform interpolation for modal $\mu$-calculus is discussed in depth in \refchapter{chapter:cyclic} of this volume.

While the semantic method via bisimulation quantifiers is highly versatile, its main drawback is that it does not give an explicit construction of the uniform interpolants of a formula. In general, the number $n$ in the Expansion Lemma (Lemma~\ref{lem:combi}) is taken in the order of several exponential factors larger than $k$, so that, while it is theoretically computable, it does not yield an efficient way for calculating uniform interpolants. We come back to this point in Section~\ref{sec:connections}.%

\section{Uniform Interpolants via Adjunction and Topology}
\label{sec:algebra}
In this section, we will show how the existence of uniform interpolants relates to the existence of adjoints between finitely generated free algebras (Subsection~\ref{subsec:adjoints}). We will again focus first on the case of intuitionistic propositional logic, and then show how to generalize this to other logics. Next, in the case of intuitionistic logic, we give a topological argument for why the required adjoints exist (Subsection~\ref{subsec:topology}). Finally, we relate this to questions in categorical logic that were at the origin of Pitts' investigations (Subsection~\ref{subsec:catlog}).

\subsection{Uniform Interpolants as Adjoints}\label{subsec:adjoints}
Recall that intuitionistic propositional logic is captured algebraically by the class of \emph{Heyting algebras}, i.e., structures $(A, \wedge, \vee, \To, \bot, \top)$ which are bounded distributive lattices having the further property that, for any $a,b \in A$, the element $a \To b$ is the maximum of the set $\{c \in A \ | \ a \wedge c \leq b\}$.  If $\p$ is any finite sequence of propositional variables, then, for any propositional formulas $\phi(\p)$ and $\psi(\p)$, their equivalence classes are elements of the \emph{free} Heyting algebra $\bH(\p)$, which we still denote by $\phi$ and $\psi$ by a slight abuse of notation. The basic connection between Heyting algebras and intuitionistic logic is the fact that an intuitionistic entailment $\phi \vdash \psi$ holds if, and only if, $\phi \leq_{\bH(\p)} \psi$.

Let $p, \q$ be a finite sequence of variables and let $\phi(p, \q)$ be a formula in these variables, which we view as an element of $\bH(p,\q)$. Note that the free algebra $\bH(\q)$ can be seen as a \emph{subalgebra} of $\bH(p,\q)$; formally, we write 
\begin{equation}\label{eq:inclusion}
i \colon \bH(\q) \hookrightarrow \bH(p,\q) 
\end{equation}
for the injective homomorphism that sends any equivalence class of formulas $\psi(\q)$ to itself, now viewed as an element of $\bH(p,\q)$. In these terms, the uniform interpolant $\rE_p(\phi)$ is an element of $\bH(\q)$ with the property that
\[ \rE_p(\phi) = \min \{ \psi \in \bH(\q) \ | \ \phi \leq_{\bH(p,\q)} i(\psi)\}.\] 
This equation says precisely that the function $\rE_p \colon \bH(p,\q) \to \bH(\q)$ is a \emph{lower adjoint} to the homomorphism $i$ of (\ref{eq:inclusion}).\footnote{Note that lower adjoints are sometimes referred to as \emph{left} adjoints, but we prefer `lower' in this context, in order to avoid confusion with the left-right terminology for uniform interpolants, which uses the \emph{opposite} convention. Below, we will move to homomorphisms between \emph{congruence lattices} rather than between the algebras themselves, and another order flip will take place, which makes the left-right terminology for uniform interpolants coherent with the one for adjoints.}

The previous paragraph shows that the existence of right uniform interpolants implies the existence of a lower adjoint to the homomorphism $i$ of (\ref{eq:inclusion}). Similarly, the existence of left uniform interpolants implies the existence of an \emph{upper adjoint} to the homomorphism $i$. A crucial idea underlying all semantic proof methods for implication-based uniform interpolation is that, conversely, if one wants to construct uniform interpolants, one may first try to construct such lower and upper adjoints, which will yield candidate right and left uniform interpolants, respectively. As already mentioned after the proof of Proposition~\ref{prop:definability} above, bisimulation quantifiers may be seen as semantic realizations of these adjoints. The following argument is thus, in a sense, an `algebraic' version of the proof of Proposition~\ref{prop:definability}.
\begin{proposition}\label{prop:adjoint}
    Let $p, \q$ be a finite set of variables. Consider the homomorphism $i$ of (\ref{eq:inclusion}).
    \begin{enumerate}
        \item The function $i$ has a lower adjoint if, and only if, every formula $\phi(p,\q)$ has a right uniform interpolant with respect to $p$.
        \item The function $i$ has an upper adjoint if, and only if, every formula $\phi(p,\q)$ has a left uniform interpolant with respect to $p$.
    \end{enumerate}
\end{proposition}
\begin{proof}
(1) The sufficiency of the condition was explained above. For the necessity, let $i^{\flat}$ be a lower adjoint for $i$ and let $\phi(p,\q)$ be arbitrary. Pick a formula $\psi(\q)$ in the equivalence class $i^{\flat}(\phi)$. By adjunction, $\phi \leq_{\bH(p,\q)} i(\psi)$, so that $\phi \vdash \psi$. Let $\theta$ be any $p$-free formula such that $\phi \vdash \theta$. First, using the (non-uniform) Craig interpolation theorem of intuitionistic propositional logic, pick an interpolating formula $\theta'$ using only variables in $\q$, such that
{\it (a)} $\phi \vdash \theta'$ and {\it (b)} $\theta' \vdash \theta$.
Property {\it (a)} implies that $\phi \leq_{\bH(p,\q)} i(\theta')$, so, by definition of lower adjoints, $i^{\flat}(\phi) \leq_{\bH(\q)} \theta'$. The choice of $\psi$ then gives that $\psi \vdash \theta'$, and now property {\it (b)}, together with the transitivity of the relation $\vdash$, yield $\psi \vdash \theta$. (2) Similar.
\end{proof}

For other logics and their corresponding varieties of algebras, analogous arguments still hold, although one needs to perform them at the level of the \emph{lattices of compact congruences} of free finitely generated algebras. In category-theoretic terms, compact congruences on finitely presented algebras may be seen as \emph{regular epimorphisms} in the category of finitely presented algebras, and this was the original approach taken in~\cite{GhiZaw2002}. An algebraic rendering of this idea, using the universal-algebraic notion of \emph{coherence}, was developed in~\cite{GooMetTsi2017, KowMet2019}. We now briefly recall how this works, referring to \refchapter{chapter:algebra} of this volume for more details.

Let $\mc{V}$ be a variety of algebras, in the sense of universal algebra, see, e.g., Section~2 of \refchapter{chapter:algebra}. We write $\mbf{F}_{\mc{V}}(\x)$ for the free $\mc{V}$-algebra over a set of variables $\x$. Recall that, in Section~3 of that chapter,  a consequence relation $\vdash_{\mc{V}}$ was associated to $\mc{V}$, and it was shown in Lemma~9 of the same chapter that, for any set of equations $\Sigma \cup \{E\}$ with variables in a set $\x$, we have $\Sigma \vdash_{\mc{V}} E$ if, and only if, the equation $E$ belongs to the congruence on $\mbf{F}_{\mc{V}}(\x)$ generated by $\Sigma$.  For any homomorphism $h \colon A \to B$ between algebras in $\mc{V}$, we have an adjoint pair of maps 
\[h^* \colon \Con A \leftrightarrows \Con B \colon h^{-1} ,\]
where $h^{-1}$ sends a congruence $\beta$ on $B$ to its inverse image $h^{-1}(\beta) := \{(a,a') \in A^2 \ | \ (ha,ha') \in \beta\}$, and $h^*$ sends a congruence $\alpha$ on $A$ to $\mathrm{Cg}^B(h[\alpha])$, i.e., the congruence on $B$ generated by the forward image $\{(ha,ha') \ | \ (a,a') \in \alpha\}$. Note that $h^*$ \emph{preserves finitely generated}\footnote{In some literature, in particular~\cite{GooMetTsi2017}, `finitely generated congruences' are also called `compact congruences', referring to the fact that they are the compact elements of the algebraic lattice $\Con A$.} \emph{congruences}, i.e., when $\alpha$ is a finitely generated congruence, then so is $h^*(\alpha)$. However, it is not the case in general that $h^{-1}$ preserves finitely generated congruences; for a counterexample in the variety of groups, see \cite[Example 3.6]{GooMetTsi2017}. It is shown in \cite[Prop.~3.5]{GooMetTsi2017} that right uniform interpolants with respect to $\vdash_{\mc{V}}$ exist if, and only if, for any sets of variables $\x, \y$, the inverse image $i^{-1}$ of the inclusion homomorphism 
\begin{equation}\label{eq:incl-V}
    i \colon \mbf{F}_{\mc{V}}(\y) \into \mbf{F}_{\mc{V}}(\x, \y) 
\end{equation}
preserves finitely generated congruences. Thus, this property in particular implies that the consequence relation $\vdash_{\mc{V}}$ has the interpolation property. It is moreover possible to `separate out' the usual interpolation property: right uniform interpolants exist if, and only if, (i) the usual interpolation property holds, and (ii) for any \emph{finite} sets of variables $\x, \y$, the inverse image $i^{-1}$ of the homomorphism $i$ of (\ref{eq:incl-V}) preserves finitely generated congruences. Writing $\KCon A$ for the semilattice of finitely generated congruences on $A$, (ii) implies that, for any finite sets of variables $\x,\y$, we have an adjunction
\[ i^* \colon \KCon \mbf{F}_{\mc{V}}(\y) \leftrightarrows \KCon \mbf{F}_{\mc{V}}(\x, \y) \colon i^{-1} \ . \] 
Note that, here, $i^{-1}$ is \emph{upper} adjoint to $i^*$, whereas in the case of Heyting algebras (Proposition~\ref{prop:adjoint}) we saw that right uniform interpolants are related to the existence of a \emph{lower} adjoint $i^{\flat}$ for the homomorphism $i$. The reason for this order switch is that any Heyting algebra $A$ is isomorphic to the \emph{opposite} of the poset $\KCon A$ of compact congruences on $A$. %

It was shown in \cite[Prop.~3.8]{GooMetTsi2017} that (ii) implies that $h^{-1}$ preserves finitely generated congruences for \emph{any} homomorphism between finitely presented algebras. Moreover, \cite[Thm.~2.3]{KowMet2019} establishes that (ii) is also equivalent to the variety $\mc{V}$ being \emph{coherent}, meaning that any finitely generated subalgebra of a finitely presented algebra is itself finitely presented. To summarize this discussion, we have the following:
\begin{theorem}[\cite{GooMetTsi2017, KowMet2019}]\label{thm:coherence}
    For any variety $\mc{V}$, the following are equivalent:
    \begin{enumerate}
        \item for any finite sets of variables $\x$, $\y$, the inverse image $i^{-1}$ of the homomorphism $i$ of (\ref{eq:incl-V}) preserves finitely generated congruences;
        \item for any homomorphism $h \colon A \to B$ between finitely presented $\mc{V}$-algebras, the inverse image $h^{-1}$ preserves finitely generated congruences;
        \item the variety $\mc{V}$ is coherent.
    \end{enumerate}
    If, moreover, the consequence relation $\vdash_\mc{V}$ has the interpolation property, then the above properties are also equivalent to:
    \begin{enumerate}
        \setcounter{enumi}{3}
        \item for any sets of variables $\x$, $\y$, the inverse image $i^{-1}$ of the homomorphism $i$ of (\ref{eq:incl-V}) preserves finitely generated congruences;        
        \item the consequence relation $\vdash_{\mc{V}}$ has right uniform interpolants.
    \end{enumerate}
\end{theorem}
Similarly, the existence of \emph{left} uniform interpolants is equivalent to the algebraic property that, for any sets $\x$ and $\y$ of variables, the direct image map $i^{*}$ has a \emph{lower adjoint}, also see \cite[Sec.~4]{GooMetTsi2017}.

\begin{example}
The variety of MV-algebras, associated to {\L}ukaciewicz multi-valued logic, is coherent, even though implication-style uniform interpolation fails. On the other hand, the variety of modal algebras is not coherent (Remark~\ref{rem:K-UI}), even though it has implication-style uniform interpolation~\cite[Sec.~6]{Ghi1995}. The variety of lattices, which algebraizes a minimal, non-distributive logic of `and' and `or', also fails to be coherent~\cite{Sch1983}, see also~\cite[Sec.~5]{KowMet2019}.
\end{example}

\subsection{Esakia Duality and Open Mappings}\label{subsec:topology}
The point of view that uniform interpolants for intuitionistic logic are given by adjoints to inclusions of free Heyting algebras enables a \emph{topological} analysis of uniform interpolants. The general idea of using topology for studying intuitionistic logic goes back to at least M.~H.~Stone~\cite{Sto1938}. 
To motivate the idea, observe first that the collection $\mathcal{O}(X)$ of open sets of any topological space $X$ is a Heyting algebra, which is even \emph{complete}. 
The Heyting implication in $\mathcal{O}(X)$ is given, for $U$ and $V$ open sets, by letting $U \To V$ be the \emph{topological interior} of the set of points of $X$ which are either in $V$ or not in $U$. 
In general, however, \emph{free} Heyting algebras on two or more generators are not complete~\cite[Thm.~4.2]{Bel1986}, and can therefore not be of the form $\mathcal{O}(X)$ for any space $X$.

To obtain \emph{all} Heyting algebras from a topological construction, one may use a theorem due to Stone~\cite{Sto1938}, which implies that any bounded distributive lattice, and thus in particular any Heyting algebra, can be represented as the collection of open \emph{and compact} subsets of a certain topological space. Homomorphisms of Heyting algebras then correspond to certain strongly continuous functions between the spaces in the other direction. We will here use a point of view on this duality due to Priestley~\cite{Pri1970} and Esakia~\cite{Esa1974, Esa2019}, who consider subsets that are clopen (i.e., closed and open) up-sets (i.e., upward closed) with respect to a partial order that is added as additional structure to the topological space. One way to view the results of Stone, Priestley, and Esakia is that they provide a \emph{canonical} choice for an embedding of a Heyting algebra into a complete Heyting algebra of the form $\mathcal{O}(X)$. Indeed, the space $X$ representing a Heyting algebra $H$ is such that $\mathcal{O}(X)$ is the \emph{ideal completion} of $H$. 

Esakia \cite{Esa1974} proved that the category of Heyting algebras is dually equivalent to a subcategory of the category of compact ordered topological spaces, in which the morphisms are the so-called \emph{continuous bounded maps}. We write $\bE(\p)$ for the Esakia space dual to the free Heyting algebra $\bH(\p)$ over set of generators $\p$, which is also known as the \emph{canonical model} for intuitionistic propositional logic. Concretely, the points of $\bE(\p)$ are the prime theories of intuitionistic propositional logic  in variables $\p$, where a \emph{theory} is a set of formulas closed under entailment and finite conjunction, and a theory $T$ is \emph{prime} if $\bot \not\in T$ and $\phi \vee \psi \in T$ implies $\phi \in T$ or $\psi \in T$. The partial order is given by inclusion of prime theories, and the topology on $\bE(\p)$ is generated by declaring that, for any $\phi \in \bH(\p)$, the set $\widehat{\phi}$ of prime theories containing $\phi$ is both closed and open (clopen). An important consequence of Esakia duality is that the assignment $\phi \mapsto \widehat{\phi}$ is a Heyting algebra isomorphism between $\bH(\p)$ and the Heyting algebra of \emph{clopen up-sets}  of the Esakia space $X$.

We call an Esakia space $X$ \emph{finitely copresented} if $X$ is order-homeomorphic to a clopen up-set of $\bE(\p)$, for some finite $\p$; this name comes from the fact that it is equivalent to saying that the Heyting algebra dual to $X$ is finitely presented. The following theorem is proved in~\cite[Thm.~2]{GooReg2017}, using methods similar to the Expansion Lemma (Lemma~\ref{lem:combi}) stated above.
\begin{theorem}\label{thm:openmapping}
    Every continuous bounded map between finitely copresented Esakia spaces is open.
\end{theorem}
From Theorem~\ref{thm:openmapping}, one can deduce in particular that the inclusion homomorphism $i \colon \bH(\q) \into \bH(p,\q)$ of (\ref{eq:inclusion}) has both a lower and an upper adjoint, and thus that uniform interpolants exist in intuitionistic logic, in light of Proposition~\ref{prop:adjoint}. We explain briefly how this works for the lower adjoint; see~\cite[Thm.~3]{GooReg2017} for more details. By Esakia duality, the homomorphism $i$ may equivalently be viewed as $p^{-1}$, where $p \colon \bE(p,\q) \to \bE(\q)$ is the continuous bounded map dual to $i$. By Theorem~\ref{thm:openmapping}, this map is open. In particular, for any formula $\phi(p,\q)$, the direct image $p[\widehat{\phi}]$ is also a clopen up-set of $\bE(\q)$, so that there must exist some formula $\epsilon$ with $\widehat{\epsilon} = p[\widehat{\phi}]$. Defining $i^{\flat}(\phi) := \epsilon$ now yields a lower adjoint for $i$. 

We remark that Theorem~\ref{thm:openmapping} was recently generalized by J. Marquès~\cite[Thm~1.4.26]{Mar2023PhD} to the class of \emph{intuitionistic compact ordered spaces}, also called \emph{generalized Esakia spaces} in~\cite[Def.~3.1]{HofNor2014}. In light of a duality of M.~Abbadini~\cite{Abb2019}, these spaces, which form a subclass of compact ordered spaces, correspond to a logic where truth values of basic propositions are taken in a \emph{continuous} domain, such as the unit interval, in the vein of {\L}ukaciewicz's many-valued logic, \Luklog. Marquès' result therefore corresponds to a uniform interpolation property for an intuitionistic version of an infinitary many-valued logic~\cite[Prop.~1.3.23]{Mar2023PhD}. A full syntactic analysis of this result has not yet been performed.

\subsection{Connections to Categorical Logic}\label{subsec:catlog}
The connection between uniform interpolants and adjoints is analogous to the important idea in categorical logic that quantifiers can be viewed as adjoints~\cite{Lam68, Law69}, and was already pointed out in Pitts' original paper~\cite[\S 4]{Pit1992}. Indeed, Pitts explains there that the existence of uniform interpolants for intuitionistic logic implies that \emph{every Heyting algebra can appear as the algebra of truth values of a model of second-order intuitionistic propositional logic}. In more syntactic terms, Pitts' result provides an \emph{interpretation} of intuitionistic second order quantifiers into the propositional fragment. We stress that this interpretation cannot be conservative: the entailment relation in $\mbf{IPC}$ is decidable, but the entailment relation in $\mbf{IPC}$ enriched with second-order quantifiers is not. A concrete example of a principle that holds for uniform interpolants and not for general second-order quantifiers is given in~\cite[Exa.~10]{Pit1992}.

The algebraic consequence of Pitts' theorem points to his original motivation for proving his theorem, which was the following open question:
\begin{question}[Pitts] \label{pittsq}
    Is every Heyting algebra the algebra of global truth values of some elementary topos?
\end{question}
An \emph{elementary topos} is a categorical model of higher order intuitionistic logic. The notion may be viewed as a generalization to the categorical level of Heyting algebras, which give an algebraic model of \emph{propositional} intuitionistic logic. A positive answer to Question~\ref{pittsq} would mean that any propositional intuitionistic theory can appear as the zeroth level of some higher order intuitionistic theory, and would thus in particular yield a new way of constructing such higher order theories. A classical construction of P.~Freyd shows that the answer to the restriction of Question~\ref{pittsq} to Boolean algebras is affirmative, see~\cite[Ex.~9.11]{Joh1977}, \cite[Rem.~4.8 and 4.9]{Pit2002}~and~\cite{Kuznetsov2024}.

In the rest of this subsection, we will show in more detail how Question~\ref{pittsq} and Theorem~\ref{thm:pitts} are related; while this has clearly been known since~\cite{Pit1992}, we are not aware of any reference where the relationship is fully explained. In order to do so, let us first quickly recall the categorical notions that play a role here; for a more pedagogical exposition of these notions from a logical perspective, we recommend, e.g.,~\cite{Str2004}.

Let $\mc{E}$ be a category. A \emph{subobject} of an object $A$ of $\mc{E}$ is an equivalence class of monomorphisms $m \colon S \rightarrowtail A$, where two monomorphisms are equivalent if each can be factored through the other. The collection $\Sub_{\mc{E}}(A)$ of subobjects of an object $A$ in $\mc{E}$ is then a partially ordered set under the factorization ordering. The association $A \mapsto \Sub_{\mc{E}}(A)$ is in fact part of a \emph{contravariant functor}: any morphism $f \colon A \to A'$ in $\mc{E}$ naturally gives an order preserving function $f^{-1} \colon \Sub_{\mc{E}}(A') \to \Sub_{\mc{E}}(A)$, by precomposition with $f$.

A category $\mc{E}$ is called an \emph{elementary topos} if \emph{(a)} $\mc{E}$ has finite limits and, \emph{(b)} for every object $A$ in $\mc{E}$, there exists a \emph{power object} $PA$ in $\mc{E}$ representing the functor sending $B$ to $\Sub_{\mc{E}}(A \times B)$; concretely, this means that there is a natural bijection between morphisms $B \to PA$ and subobjects of $A \times B$.
Note that, taking $A$ to be the terminal object $1$ of $\mc{E}$, the object $\Omega := P1$ is then a \emph{subobject classifier} in $\mc{E}$, in the sense  that, for any object $B$ in $\mc{E}$, morphisms from $B$ to $\Omega$ are in natural bijection with subobjects of $B$. For $S$ a subobject of $B$, we denote the corresponding morphism by $\chi_B \colon B \to \Omega$, and we call it the \emph{characteristic map} of $B$. The prototypical example of an elementary topos is the category of sets, where finite limits can be computed as subsets of Cartesian products, $PA$ is the power set of $A$, and $\Omega$ is the two-element set.

The two axioms for elementary toposes have many consequences, of which we mention two important ones here. First, for any object $A$ in an elementary topos $\mc{E}$, the poset $\Sub_{\mc{E}}(A)$ is a Heyting algebra, and for any morphism $f \colon A \to A'$ in $\mc{E}$, the morphism $f^{-1}$ is a Heyting algebra homomorphism, see, e.g.,~\cite[Thm.~13.2]{Str2004}.
In particular, taking the terminal object $1$ of $\mc{E}$, we have a Heyting algebra $\Sub_{\mc{E}}(1)$, which is called the \emph{algebra of global truth values}.
Question~\ref{pittsq} now asks if any Heyting algebra may appear as such an algebra of global truth values. One may show that, in order to answer this question in the positive, it would in fact suffice to establish this for free Heyting algebras~\cite[p.~37]{Pitts2025}; we omit the proof.

We will now explain, in Proposition~\ref{prop:eltop-quants}, why a positive answer to 
Question~\ref{pittsq} would immediately yield Theorem~\ref{thm:pitts}, as stated in~\cite[p.~36]{Pit1992}.
To this end, recall that, in an elementary topos $\mc{E}$, there exist existential and universal quantifiers with respect to any morphism. More precisely, for any morphism $f \colon A \to B$ in $\mc{E}$, the morphism $f^{-1} \colon \Sub_{\mc{E}}(B) \to \Sub_{\mc{E}}(A)$ has an upper adjoint $\forall_f$ and a lower adjoint $\exists_f$, see, e.g.,~\cite[Thm.~13.4 and Thm.~13.5.3]{Str2004}. In the proof of the next proposition, we will apply this fact to the (unique) morphism $! \colon \Omega \to 1$.
\begin{proposition}\label{prop:eltop-quants}
    Let $\q$ be an infinite set of variables and $p \not\in \q$. Suppose that $\mc{E}$ is an elementary topos with an isomorphism $\tau \colon \Sub_{\mc{E}}(1) \cong \mbf{H}(\q) \cocolon \sigma$. 
    There exists a Heyting algebra embedding $j \colon \mbf{H}(p,\q) \hookrightarrow \Sub_{\mc{E}}(\Omega)$ such that, for any formula $\phi(p,\q)$, the formulas $\tau(\exists_{!}(j(\phi)))$ and $\tau(\forall_{!}(j(\phi)))$ are right and left uniform interpolants of $\phi$ with respect to $p$.
\end{proposition}
\begin{proof}
    For each $q \in \q$, let $j(q)$ be the subobject $!^{-1}(\sigma(q))$ of $\Omega$, and moreover let $j(p)$ be the subobject of $\Omega$ corresponding to the identity morphism $\Omega \to \Omega$. These assignments extend uniquely to a homomorphism $j \colon \mbf{H}(p,\q) \to \Sub_{\mc{E}}(\Omega)$.
    We first show that $j$ is injective.
    For this, consider the function $\alpha \colon \Sub_{\mc{E}}(\Omega) \to \mbf{H}(\q)^{\mbf{H}(\q)}$, defined, for $S \in \Sub_{\mc{E}}(\Omega)$, by letting $\alpha(S) \colon \mbf{H}(\q) \to \mbf{H}(\q)$ be the function that sends $\phi \in \mbf{H}(\q)$ to 
    $\tau(\chi_{\sigma(\phi)}^{-1}(S))$.
    Observe that this function $\alpha$ is a homomorphism of Heyting algebras. %
    One may moreover verify that, for each $q \in \q$, $\alpha(j(q))$ is the constant function with value $q$, and that $\alpha(j(p))$ is the identity function. 
    Using that $\alpha$ and $j$ are homomorphisms of Heyting algebras, it follows that, for any $\psi \in \mbf{H}(p,\q)$, $\alpha(j(\psi))$ is the function which sends $\phi \in \mbf{H}(\q)$ to the formula $\psi(\phi/p, \q)$, in which $\phi$ is substituted for all occurrences of $p$. Now, let $\psi,\psi' \in \mbf{H}(p,\q)$ be arbitrary, and suppose that $j(\psi) = j(\psi')$. Pick a variable $r$ from $\q$ which does not appear in $\psi$ or $\psi'$. Then $\psi(r/p, \q) = \alpha(j(\psi)) = \alpha(j(\psi')) = \psi'(r/p,\q)$, and thus $\psi = \psi'$.

    To finish the proof, consider the commutative diagram
\begin{center}
    \begin{tikzcd}
        \mbf{H}(\q) \arrow[r, "i"] \arrow[d, "\sigma"] & \mbf{H}(p, \q) \arrow[d, "j"]\\
        \Sub_{\mc{E}}(1) \arrow[r, "!^{-1}"] & \Sub_{\mc{E}}(\Omega)
    \end{tikzcd}
\end{center}
    Note that $\tau \circ \exists_!$ is lower adjoint to $!^{-1} \circ \sigma = j \circ i$. Thus, since $j$ is an order embedding, we get that $\tau \circ \exists_{!} \circ j$ is lower adjoint to $i$, and we conclude by Proposition~\ref{prop:adjoint}.
\end{proof}

In an attempt to settle Question~\ref{pittsq} in the negative, in light of Proposition~\ref{prop:eltop-quants}, Pitts had initially tried to prove the negation of Theorem~\ref{thm:pitts}. However, Theorem~\ref{thm:pitts} ended up being true, and Question~\ref{pittsq} thus remained open. Theorem~\ref{thm:pitts} has since been used in attempts towards settling Question~\ref{pittsq} in the positive, in particular by the late D.~Pataraia, although there is no published proof; see, e.g.,~\cite{Joh2012, MamukaSO2017, Kuznetsov2024, Pitts2025} for more information.

The above category-theoretic viewpoint on uniform interpolation can also be related to the model-theoretic results that can be derived from it, as we will describe in the next section.

\section{Model-complete Theories from Uniform Interpolants}
\label{sec:model-theory}

In the algebraic viewpoint on propositional logics, we typically study a logical entailment relation as a special case of the equational consequence relation $\vdash_{\mc C}$ for a class of algebras $\mc C$ in an algebraic language $\mc{L}_{\mc{C}}$. For instance, the entailment relation of intuitionistic propositional logic can be viewed in this way as $\vdash_{\mc{H}}$, where $\mc{H}$ is the class of Heyting algebras, which are structures in the algebraic language ${\mc L}_{\mc H} = \{\wedge, \vee, \To, \bot, \top\}$.

In this section, rather than looking at $\vdash_{\mc{C}}$ directly, we take a different perspective, and study the \emph{first order theory of} a class of algebras $\mc C$ itself, which is defined as
\[  \ThFO(\mc{C}) := \{ \Phi \text{ a first order $\mc{L}_{\mc{C}}$-sentence such that } A \models \Phi \text{ for every } A \in \mc{C} \}.  \] 
For example, $\ThFO(\mc H)$ contains the sentences $\forall x. \forall y. x \vee y = y \vee x$ and $\exists y. \forall x. x \wedge y = y$, but it does not contain the sentence $\forall x. \exists y. ((x \vee y = \top) \mand (x \wedge y = \bot))$, as witnessed, for instance, by the three-element chain.

\begin{notation}
    Note that we write $\mand$ for the conjunction of first order logic,  in order to avoid overloading the symbol $\wedge$, which is already used as an operation symbol of the language $\mathcal{L}_{\mathcal{H}}$ of Heyting algebras; we adopt similar conventions for the other connectives of first order logic. Moreover, in this section we will use the word \emph{term} to refer to a propositional formula, since it is exactly the same thing as a term in the language $\mc{L}_{\mc{H}}$, while we reserve the word \emph{formula} for formulas of first order logic; as usual, a \emph{sentence} is a first order formula in which all the variables are bound. In order to clearly separate metavariables for propositional formulas from those for first order formulas, as a convention, we continue to use lower case Greek letters to refer to propositional formulas, while we reserve upper case Greek letters for formulas of first order logic.
\end{notation}

The starting point of the work of S. Ghilardi and M. Zawadowski in~\cite{GhiZaw1997, GhiZaw2002} is that Pitts' Theorem~\ref{thm:pitts} can be used to establish that \emph{the first order theory of Heyting algebras has a model completion}~\cite[p.~29]{GhiZaw1997}. The concept of model completion that is used here originates in A.~Robinson's model-theoretic algebra~\cite{Rob1951, Rob1963}, and closely relates to the idea of \emph{quantifier elimination}.
The prime example is the theory of algebraically closed fields, which is the model completion of the first order theory of integral domains. An example closer to the kind we want to consider here is the first order theory of Boolean algebras, which has a model completion, namely, the theory of {atomless} Boolean algebras, see, e.g.,~\cite[Thm.~6.21]{Poi2000}. 
We now recall the definitions and some basic facts about these concepts, referring to, e.g.~\cite[Sec.~3.5]{ChaKei1990} for more details. %

For the remainder of this section, we fix an algebraic language $\mc{L}$. A (first order) \emph{theory} is a set $T$ of first order sentences in the language $\mc{L}$, and a \emph{model} of $T$ is an $\mc{L}$-structure $A$ such that $A \models \Phi$ for every $\Phi \in T$. When $T$ is a theory and $\Phi(\x)$ is a formula,  we also write $T \vdash \Phi$ if, for any structure $A$ and any $\a \in A^{\x}$, if $A \models \Psi(\x/\a)$ for every $\Psi \in T$, then $A \models \Phi(\x/\a)$. 

We call a formula $\Phi$ \emph{universal} if its prenex form only contains universal quantifiers and we similarly define the notion of \emph{existential formula}. We say that two theories $T$ and $T'$ in the same language are \emph{co-theories} if, for any universal sentence $\Phi$, we have $T \vdash \Phi$ if, and only if, $T' \vdash \Phi$. In other words, writing $T_{\forall}$ for the set of universal sentences $\Phi$ such that $T \vdash \Phi$, $T$ and $T'$ are co-theories means that $T_{\forall} = T'_{\forall}$. A theory $T$ is \emph{universally axiomatizable}, or simply \emph{universal}, if, for any $\Phi \in T$, we have $T_{\forall} \vdash \Phi$. 
We say that a theory $T$ is \emph{model complete} if, for any formula $\Phi$, there is a universal formula $\Phi'$ such that $T \vdash \Phi \leftrightarrow \Phi'$. 
Notice that it is equivalent to say that, for any formula $\Psi$, there is an existential formula $\Psi'$ such that $T \vdash \Psi \leftrightarrow \Psi'$,  as can be seen by applying the definition of model-completeness to $\neg \Psi$. %
A stronger property than model-completeness is that the theory $T$ \emph{eliminates quantifiers}, which means that,  for any formula $\Phi$, there is a quantifier-free formula $\Phi'$ such that $T \vdash \Phi \leftrightarrow \Phi'$. 

\begin{mydefinition}
    Let $T, T_0$ be theories in the language $\mc{L}$, and suppose that $T_0$ is universal.\footnote{There is a more general definition where $T_0$ is not required to be universal, but it is slightly more complicated to state, and not relevant for the examples of interest to us here.} 
    \begin{itemize}
        \item The theory $T$ is a \emph{model companion} of $T_0$ if $T$ is a model complete co-theory of $T_0$. 
        \item The theory $T$ is a \emph{model completion} of $T_0$ if $T$ is a co-theory of $T_0$ and eliminates quantifiers.
    \end{itemize}
\end{mydefinition}

There is a complementary semantic perspective on model completions, which we also briefly recall now. Let $T_0$ be a universal theory and let $A$ be a model of $T_0$. The model $A$ is called \emph{existentially closed} for $T_0$ provided that, for any quantifier free formula $\Phi(\x,y)$ and $\a \in A^{\x}$, if there exists a model $B$ of $T_0$ such that $B$ contains $A$ as a submodel and $B \models \exists y. \Phi(\x/\a, y)$, then $A \models \exists y. \Phi(\x/\a,y)$. Note, for instance, that the existentially closed integral domains are precisely the algebraically closed fields in the usual sense. In general, if $T$ is a model companion of a universal theory $T_0$, then the models of $T$ are precisely the {existentially closed} models of $T_0$. From this, one can deduce that the universal theory $T_0$ has at most one model companion, up to equivalence, and that this model companion exists if, and only if, there is a first order theory whose models are exactly the existentially closed models of $T_0$. 
Moreover, one may show that this model companion is a model completion if, and only if, the class of models of $T_0$ has the amalgamation property.

We now sketch how a model completion for the theory of Heyting algebras may be obtained from Pitts' quantifiers~\cite{GhiZaw1997}.

\begin{theorem}\label{thm:HA-model-completion}
    The first order theory of Heyting algebras, $\ThFO(\mc H)$, has a model completion.
\end{theorem}
\begin{proof}[Proof (Sketch).]
    It is well-known that the class of Heyting algebras has the amalgamation property, see for instance~\cite[Prop.~1]{Mak1977} or~\cite{Fle1980} for early algebraic proofs of this fact. In light of the above remarks, it therefore suffices to show that there is a first order theory whose models are exactly the existentially closed Heyting algebras.

Let $\x, y$ be a finite set of variables  and let $\phi_1(\x,y), \dots, \phi_n(\x,y)$ and $\psi_1(\x,y), \dots, \psi_m(\x,y)$ be any terms in these variables, and consider the quantifier-free formula $\Psi(\x)$ given by  
\[ \exists y. (\phi_1 = \top \mand \cdots \mand \phi_n = \top \mand \psi_1 \neq \top \mand \cdots \mand \psi_m \neq \top)\ . \] 
We will show how to express in a first order way that a Heyting algebra satisfies the existential closedness property with respect to this formula $\Psi$; the case of a general existential formula is then essentially an iteration of this argument.
    Consider the formula $\Psi'(\x)$ given by
\[  \rE_y\left(\bigwedge_{i=1}^n \phi_i\right) = \top \mand \rA_y\left( \left[\bigwedge_{i=1}^n \phi_i\right] \Rightarrow \psi_1\right) \neq \top \mand \cdots \mand  \rA_y \left(\left[\bigwedge_{i=1}^n \phi_i\right] \Rightarrow \psi_m \right) \neq \top\ , \]
where we use the uniform interpolants for propositional formulas, which exist according to Theorem~\ref{thm:pitts}.
Using the definition of uniform interpolants, one may then prove that, for any Heyting algebra $A$ and $\a \in A^{\x}$, we have $A \models \Psi'(\x/\a)$ if, and only if, 
$A[y]/\alpha \models \Psi(\x/\a)$,
where $A[y]$ denotes the coproduct of $A$ with the free Heyting algebra over $y$, and $\alpha$ is the congruence generated by the pairs $(\phi_i(\x/\a,y), \top)$, for $i = 1, \dots, n$.
            Therefore, a Heyting algebra $A$ satisfies the existential closedness property for the formula $\Psi(\x,y)$ and a tuple $\a \in A^{\x}$ if, and only if, $A \models (\Psi' \rightarrow \Psi)(\x/\a)$. Thus, an existentially closed Heyting algebra must satisfy the first order sentence $\forall \x. (\Psi' \rightarrow \Psi)$. 
            This is the essential property required for axiomatizing existentially closed Heyting algebras. 
            Using this property, a general argument 
            (see, e.g., \cite[Thm.~3.5.20]{ChaKei1990}) 
            shows that the model completion of $\mathrm{Th}_{\mathsf{FO}}(\mathcal{H})$ is obtained by adding to the theory $\mathrm{Th}_{\mathsf{FO}}(\mathcal{H})$, for any choice of $\x, y, \phi_1, \dots, \phi_n, \psi_1, \dots, \psi_m$ as above,  the formula $\forall \x. (\Psi' \to \Psi)$.
\end{proof}%
The above proof method shows that the class of existentially closed Heyting algebras admits a first order axiomatization $T$, and also gives an effective quantifier elimination procedure for $T$.  \cite[Prop.~A.2]{GhiZaw1997} establishes some concrete properties of existentially closed Heyting algebras, such as the fact that the order must be dense, and that there are no meet-irreducible elements. It is an open problem to obtain a `manageable' complete axiomatization of the class of existentially closed Heyting algebras.  \cite[Thm.~A.3]{GhiZaw1997} further shows that all of the eight varieties of Heyting algebras with the amalgamation property also have model completions. L.~Darnière and M.~Junker~\cite{DarJun2018} obtain concrete axiomatizations for the model completions of all of these varieties, \emph{except} for the variety of all Heyting algebras. 
Several more results on the (non-)existence of model completions for classes of modal algebras are obtained in \cite{GhiZaw2002}. More generally, an algebraic method for establishing (non-)existence of model companions for many modal logics, substructural logics, and varieties of lattices, is developed in~\cite{KowMet2018,KowMet2019}.

We remark that Theorem~\ref{thm:HA-model-completion} can be viewed as a consequence of a more general phenomenon that goes back to work of W. Wheeler~\cite{Whe1976, Whe1978}, as explained in~\cite{GhiZaw1997,GhiZaw2002}, and also further investigated in~\cite{GooMetTsi2017, MetReg2023}. The general idea is that, for $T$ a universal theory in an algebraic language $\mc{L}$, the existence of a model companion for $T$ can be characterized via properties of the category $\mc{F}_{T}$ of \emph{finitely presentable $T$-structures}, with homomorphisms between them. For instance, it is proved in~\cite[Thm.~3.4]{GhiZaw1997} that, for an equational theory $T$, if the opposite $(\mc{F}_T)^{\op}$ is an \emph{r-Heyting category}, then $T$ has a model completion. Moreover, in the presence of a few additional assumptions on $T$, the converse also holds. From this perspective, a generalized form of Theorem~\ref{thm:pitts} can first be used to establish that the opposite of the category of finitely presented Heyting algebras is an r-Heyting category, from which Theorem~\ref{thm:HA-model-completion} then follows by the general theory; this is the approach taken in \cite{GhiZaw2002}. Using the converse direction, one may also use this general approach to deduce that certain theories \emph{do not} have a model completion, as is for instance done for many extensions of modal logic $\mbf{K4}$ in \cite[Ch.~6]{GhiZaw2002}.

\section{Syntactic Methods, Further Results and Connections}
\label{sec:connections}
Our focus in this chapter so far has been on connections between uniform interpolation and \emph{semantics}. The \emph{syntactic} method of establishing uniform interpolation is discussed in more detail in Section~3 of~\refchapter{chapter:prooftheory} in this volume; we will here just make some brief remarks about the method and its relation to the semantic method. 

The basic idea of the syntactic proof style is to recursively construct candidate formulas $\rE_p(\phi)$ and $\rA_p(\phi)$ for a formula $\phi$, usually based on a large case distinction on the shape of $\phi$, and the rules that may be applied to it in a particular sequent calculus for the logic. The most intricate part of such proofs of uniform interpolation is the proof of correctness of these candidate formulas $\rE_p(\phi)$ and $\rA_p(\phi)$, which is typically done by induction on the structure of a `potential derivation' using $\phi$ in the sequent calculus. 
In order to make such an induction possible, one needs a well-founded order on the set of potential derivations that use $\phi$.
For many logics, including $\mathbf{IPC}$, it is not a priori obvious how to obtain a sequent calculus that terminates in the strong sense required for constructing uniform interpolants. 

The syntactic method was in particular the one that was originally used by Pitts to establish Theorem~\ref{thm:pitts} in~\cite{Pit1992}, and was subsequently used for proving the implication-based uniform interpolation property for a great number of other logics; for instance, this was done by Bílková~\cite{Bil2007} for the basic modal logic $\mathbf{K}$ and for the modal logic of reflexive frames, $\mathbf{T}$, and for provability logics~$\mathbf{GL}$ and $\mathbf{Grz}$ in~\cite{BilPhD}, also see~\cite{Bil22}. The method was further extended by Iemhoff~\cite{Iem2019ARCH, Iem2019APAL}, who identifies precise general criteria for sequent calculi that allow one to establish uniform interpolation; these ideas were subsequently applied to obtain syntactic proofs of uniform interpolation for a number of intuitionistic modal logics, see, e.g.,~\cite{Gie22, IemJalTab2022, ShiGieGorIem2023}, and the further references given in \refchapter{chapter:prooftheory} of this volume. In the context of substructural logic, Alizadeh, Derakhshan, and Ono~\cite{AliDerOno2014} establish uniform interpolation for full Lambek calculus and for linear and affine logics without exponentials. Interestingly, in that context, even the \emph{predicate} versions of the logics have uniform interpolation, in contrast with the case of classical logic, see e.g. Theorem~2.17 in \refchapter{chapter:firstorder}.

As mentioned at the end of Section~\ref{sec:bisimulation}, an advantage of this syntactic method over the semantic one is that, at least in theory, it provides better bounds on the complexity of the uniform interpolants. The intuitive reason for this is that, to find a uniform interpolant for a formula $\phi$, the semantic method performs a search among `all possible bisimilarity types of Kripke models, considered up to the depth of $\phi$', while the syntactic method searches `all possible proofs of $\phi$'. The latter is in general a (much) smaller number. In practice, however, it is usually not feasible to compute uniform interpolants by hand in a logic such as $\mbf{IPC}$, even with the syntactic method, as the calculations quickly become complex even on small examples. 

Uniform interpolation is also studied in the context of knowledge representation via description logics, where it is in particular relevant for ontology extensions~\cite{GhiLutWol2006}, and where the complexity and size of uniform interpolants have been studied in detail; see~\refchapter{chapter:kr} of this volume.

The algorithms for computing uniform interpolants are often intricate, and it is a non-trivial task to implement them correctly. To mitigate this, in~\cite{FerGoo2023}, a version of Pitts' algorithm for computing uniform interpolants for $\mbf{IPC}$ was implemented and moreover verified to be correct in the Rocq proof assistant. These methods were subsequently extended to provide a verified computation of uniform interpolants for the modal logics $\mbf{K}$ (where the uniform interpolants are not deductive but implication-based), for $\mbf{GL}$, and for the intuitionistic modal logic $\mbf{iSL}$~\cite{FGGS2024}. This work in particular provides a uniform interpolant calculator which is available as an online application at \url{https://hferee.github.io/UIML/}, and has recently been used in theoretical research for computing specific uniform interpolants~\cite{Koc2025}.

\begin{table}[t]
    \caption{Summary of the literature on uniform interpolation and model companions. The results on Boolean algebras are folklore. A double \xmark~means that even the \emph{non-uniform} form of interpolation fails in the logic. An equality sign in the `deductive UI' column means that a deduction property holds in the logic, so that 
the result is the same as for implication-based UI; see the third remark following Definition~\ref{dfn:UI}. In the last column, the symbol ${}^{A}$ indicates that the class of algebras also has amalgamation, so that the model companion is in fact a model completion. Unless mentioned otherwise, the results in the last column are due to~\cite{GhiZaw1997}.}
    \label{tab:summary}
  \centering
    \begin{tabularx}{\textwidth}{X|X|X|X}
        logic & {implication UI} & deductive UI & {FO-theory has model companion}\\ \hline \hline
        $\mathbf{CPC}$ & {\cmark} & = & \cmark${}^{A}$ [folklore] \\ \hline
        $\mathbf{IPC}$ & \cmark~\cite{Pit1992} & = &  \cmark${}^{A}$ \\ \hline
        $\mathbf{K}$ & \cmark~\cite[Sec.~6]{Ghi1995} & \xmark~\cite{GhiLutWol2006long}, cf. Remark~\ref{rem:K-UI} &  \xmark \\ \hline
        $\mathbf{S4}$ & \xmark~\cite{GhiZaw1995a} & = &  \xmark~\cite[Thm.~2]{Lip1982} \\ \hline
        $\mathbf{K4}$ & \xmark~\cite{Bil2007} & = & \xmark \\ \hline 
        $\mathbf{GL}$ & \cmark~\cite[Cor.~2.13]{Shavrukov} & = & \cmark${}^{A}$ \\ \hline
        $\mathbf{Grz}$ & \cmark~\cite[Sec.~8]{Vis1996} & = &  \cmark${}^{A}$ \\ \hline
        {\Luklog}  & \xmark\xmark~\cite{KrzZac1977, Mun2011a} & \cmark~\cite[Rem.~4]{Mun2011a} & \xmark~\cite{Lac1979} \\ \hline
        $\mathbf{LTL}$ & \xmark\xmark~\cite{Mak1991} & = & \cmark~(a variant)~\cite{GhiGoo2016}\\ \hline
        $\mathbf{K}_{\mu}$ & = & \cmark~\cite{AgoHol2000} &  \\  \hline
        $\mathbf{iK, iKD}$ & \cmark~\cite{Iem2019APAL} &   & \\  \hline
        $\mathbf{iSL}$ & \cmark~\cite{LitakVisser24, FGGS2024} & &  \\ \hline
        $\mathbf{iK4}, \mathbf{iS4}$ & \xmark~\cite{Gie22} & &  \\  \hline
        $\mathbf{iGL}, \mathbf{IK}, \mathbf{CK}$ &  & &  \\ \hline
    \end{tabularx}
\end{table}

In Table~\ref{tab:summary}, we summarize the current state of the art of the literature on uniform interpolation for a number of common logics, and the existence of model companions for their corresponding classes of algebras. Note that this table also implies a number of open questions: for all the empty cells, at the time of writing, the author is unaware of any answers, yet hopeful that readers of this chapter will be inspired to produce some.\footnote{Just before this chapter was going to press in December 2025, the author was informed that uniform interpolation for iGL has been proved by I. van der Giessen, B. Sierra Miranda, and G. Menéndez Turata.} For more information, we refer to the table maintained by T. Kurahashi on the webpage \url{https://www2.kobe-u.ac.jp/~tk/jp/notes/ULIP.html}, which surveys Craig interpolation, implication-style uniform interpolation, and their Lyndon variants, for a large number of intermediate, modal, and intuitionistic modal logics.

\section*{Acknowledgments}
The author is grateful for several discussions on parts of this chapter with Balder ten Cate, Silvio Ghilardi, Iris van der Giessen, Jérémie Marquès, George Metcalfe, and Frank Wolter. Van der Giessen and Metcalfe also each provided detailed reviews of earlier versions of this chapter, which have been a great help in improving it.
This work has received financial support from the Agence Nationale de la Recherche (ANR), project ANR-23-CE48-0012-01.
\bibliography{uichapter.bib,taci}

\end{document}